\documentclass[leqno,a4paper,11pt]{article}
\usepackage{amssymb}

\usepackage[cmex10]{amsmath}
\usepackage{amsfonts}
\usepackage{amsthm}

\usepackage{fixltx2e}

\usepackage{url}	  

\usepackage{lmodern}

\newcommand{\kr}{\mathrm{ker} \,}

\theoremstyle{definition} \newtheorem{theorem}{Theorem}[section]
\theoremstyle{definition} \newtheorem{definition}[theorem]{Definition}
\theoremstyle{definition} \newtheorem{lemma}[theorem]{Lemma}
\theoremstyle{definition} 
\theoremstyle{definition} \newtheorem{corollary}[theorem]{Corollary}
\theoremstyle{definition} \newtheorem{observation}[theorem]{Observation}
\theoremstyle{definition} \newtheorem{example}[theorem]{Example}
\theoremstyle{definition} 
\theoremstyle{definition} \newtheorem{remark}[theorem]{Remark}
\theoremstyle{definition} 
\newcommand{\eat}[1]{}

\begin{document}

\title{$(m,n)$-Semirings and a \\Generalized Fault Tolerance Algebra of Systems}			
\author{Syed Eqbal Alam~\footnote{\tt syed.eqbal@iiitb.net}
\\ College of Computers and Information Technology \\
Taif University \\
Taif, Saudi Arabia
\and Shrisha Rao~\footnote{\tt shrao@alumni.cmu.edu} 
    \\
IIIT - Bangalore \\ 
    Bangalore 560 100, India
\and Bijan Davvaz~\footnote{\tt davvaz@yazduni.ac.ir}
\\
Department of Mathematics \\
Yazd University,
\\
 Yazd, Iran
 }

\date{}

\maketitle						

\begin{abstract}

  We propose a new class of mathematical structures called
  \emph{$(m,n)$-semirings} (which generalize the usual semirings), and
  describe their basic properties.  We define partial ordering and
  generalize the concepts of congruence, homomorphism, etc., for
  $(m,n)$-semirings.  Following earlier work by
  Rao~\cite{SRaoAlgeb2008}, we consider systems made up of several
  components whose failures may cause it to fail, and represent the
  set of such systems algebraically as an $(m,n)$-semiring.  Based on
  the characteristics of these components, we present a formalism to
  compare the fault-tolerance behavior of two systems using our
  framework of a partially ordered $(m,n)$-semiring.

\end{abstract}

{\bf Keywords:} $(m,n)$-semiring, components, fault tolerance, partial
ordering
  \footnote{The work was done while first author was a Masters student at International Institute
of Information Technology(IIIT)-Bangalore}  

\section{Introduction}

Fault tolerance is the property of a system to be functional even if
some of its components fail.  It is a very critical issue in the
design of the systems as in Air Traffic Control
Systems~\cite{Cristian96,briere1993}, real-time embedded
systems~\cite{TolgaAyav2008},
robotics~\cite{Ferrell1994,Benjamin2007}, automation
systems~\cite{Avizienis1987,Cristian1990}, medical
systems~\cite{capua2007}, mission critical systems~\cite{Perraju96}
and a lot of others.  Description of fault tolerance modeling using
algebraic structures is proposed by Beckmann~\cite{beckmann1992} for
groups, and by Hadjicostis~\cite{hadjicostis1995} for semigroups and
semirings.  Semirings are also used in other areas of computer science
like cryptography~\cite{christopher2002}, databases~\cite{todd2007},
graph theory, game theory~\cite{golan1999}, etc.
Rao~\cite{SRaoAlgeb2008} uses the formalism of semirings to analyze
the fault-tolerance of a system as a function of its composition, with
a partial ordering relation between systems used to compare their
fault-tolerance behaviors.
 	 
The generalization of algebraic structures were in active research for
a long time, Timm~\cite{timm1967} in 1967 proposed commutative n-groups, later 
\newline
Crombez~\cite{Crombez1972} in 1972 generalized rings and named it as $(m,n)$-rings. It was further studied by  Crombez  and Timm~\cite{CrombezTimm1972},
Leeson and Butson~\cite{LeesonBut1980, leesonbutson1980} and by Dudek~\cite{Dudek1981}. Recently the generalization of algebraic structures is studied by Davvaz, Dudek, Vougiouklis and Mirvakili ~\cite{davvazdudekvoug2009,davvazdudekmir2009}. 

 In this paper, we first define the $(m,n)$-semiring
\((\mathcal{R},f,g)\) (which is a generalization of the ordinary
semiring \((\mathcal{R},+,\times)\), where \(\mathcal{R}\) is a set
with binary operations $+$ and $\times$), using $f$ and $g$ which are
$m$-ary and $n$-ary operations respectively.  We propose identity
elements, multiplicatively absorbing elements, idempotents, and
homomorphisms for $(m,n)$-semirings.  We also briefly touch on
zero-divisor free, zero-sum free, additively cancellative, and
multiplicatively cancellative $(m,n)$-semirings, and the congruence
relation on $(m,n)$-semirings.  In Section~\ref{partial_ordering} we
use the facts that each system consists of components or sub-systems,
and that the fault tolerance behavior of the system depends on each of
the components or sub-systems that constitute the system.  A system
may itself be a module or part of a larger system, so that its
fault-tolerance affects that of the whole system of which it is a
part.  We analyze the fault tolerance of a system given its
composition, extending earlier work of
Rao~\cite{SRaoAlgeb2008}.  Section~\ref{prelim} describes the notations
used and the general conventions followed.

Section~\ref{sem_sec} deals with the definition and properties of
$(m,n)$-semirings.  In Section~\ref{partial_ordering} we extend the
results of Rao~\cite{SRaoAlgeb2008} using a partial ordering on the
$(m,n)$-semiring of systems: class of systems is algebraically
represented by an $(m,n)$-semiring, and the fault tolerance behavior
of two systems is compared using partially ordered $(m,n)$-semiring.

\section{Preliminaries} \label{prelim}

The set of integers is denoted by $\mathbb{Z}$, with $\mathbb{Z}_+$
and $\mathbb{Z}_-$ denoting the sets of positive integers and negative
integers respectively and $m$ and $n$ used are positive integers. Let $\mathcal{R}$ be a set and $f$ be a
mapping $f:\mathcal{R}^m \to \mathcal{R}$, i.e., $f$ is an $m$-ary
operation. Elements of the set $\mathcal{R}$ are denoted by $x_i,y_i$
where $i \in \mathbb{Z}_+$.

\begin{definition} \label{n_ary_groupoid} 

  A nonempty set $\mathcal{R}$ with an $m$-ary
  operation $f$ is called an \emph{$m$-ary groupoid} and is denoted by
  ($\mathcal{R}$,$f$) (see Dudek~\cite{dudek2001}).

\end{definition}

We use the following general convention:

The sequence $x_{i}$, $x_{i+1}$, \ldots, $x_{m}$ is denoted by
$x_i^{m}$ where $1 \leq i \leq m$. 

For all $1\leq i \leq j \leq m$, the following term:
\begin{equation} \label{defeqnfm}
f(x_{1}, \ldots ,x_{i},y_{i+1}, \ldots ,y_{j},z_{j+1} , \ldots , z_{m}) 
\end{equation}
is represented as: 
\begin{equation} \label{defeqnoff}
f(x_1^{i},y_{i+1}^{j},z_{j+1}^{m}) 
\end{equation}

In the case when $y_{i+1} = \ldots = y_{j} = y$, \eqref{defeqnoff} is
expressed as:
\[f(x_{1}^i , \stackrel{(j-i)}{y} , z_{j+1}^m)\] 

\begin{definition} \label{definition_additive} 
 Let $x_1,x_2, \ldots,x_{2m-1}$ be elements of set $\mathcal{R}$. 

\begin{itemize}
\item[(i)]
  Then the associativity and
  distributivity laws for the $m$-ary operation $f$ are defined as
  follows:

\begin{itemize}

\item[(a)] \emph{Associativity}:

  $$f(x_1^{i-1}, f(x_i^{m+i-1}), x_{m+i}^{2m-1}) = f(x_1^{j-1},
  f(x_j^{m+j-1}), x_{m+j}^{2m-1}),$$

for all $x_1, \ldots, x_{2m-1} \in \mathcal{R}$, for all $1 \leq i
\leq j \leq m$ (from Gluskin~\cite{gluskin1965}).
  
\item[(b)] \emph{Commutativity}:
  $$f(x_1 , x_2 , \ldots , x_m) = f(x_{\eta(1)} , x_{\eta(2)} ,
  \ldots , x_{\eta(m)}),$$

  for every permutation $\eta$ of $\{1,2, \ldots, m\}$ (from
  Timm~\cite{timm1967}), $\forall$ $x_1$,$ x_2$,$ \ldots, x_m$ $\in
  \mathcal{R}$.
\end{itemize}

\item[(ii)] An $m$-ary groupoid $(\mathcal{R}, f)$ is called an
  \emph{$m$-ary semigroup} if $f$ is associative (from
  Dudek~\cite{dudek2001}); i.e., if
   
 $$f(x_1^{i-1}, f(x_i^{m+i-1}), x_{m+i}^{2m-1}) = f(x_1^{j-1},
  f(x_j^{m+j-1}), x_{m+j}^{2m-1}),$$

for all $x_1, \ldots, x_{2m-1} \in \mathcal{R}$, where $1 \leq i \leq
j \leq m$.

 \item[(iii)] Let $x_{1},x_2,\ldots, x_{n}, a_{1}, a_2, \ldots, a_{m}$
   be elements of set $\mathcal{R}$, and $1 \leq i \leq n$.  The
   $n$-ary operation $g$ is \emph{distributive} with respect to the
   $m$-ary operation $f$ if
  $$g(x_1^{i-1}, f(a_{1}^{m}), x_{i+1}^{n}) = f(g(x_1^{i-1},
   a_{1},x_{i+1}^{n}), \ldots, g(x_1^{i-1}, a_{m}, x_{i+1}^{n})).$$

\end{itemize}
\end{definition}

\begin{remark}
\begin{itemize}

\item[(i)]

An $m$-ary semigroup $(\mathcal{R},f)$ is called a \emph{semiabelian}
or $(1,m)$-commutative if
$$f(x, \underbrace{a, \ldots, a}_{m-2}, y) = f(y, \underbrace{a,
  \ldots, a}_{m-2}, x)$$ for all $x, y, a \in \mathcal{R}$ (from Dudek
and Mukhin~\cite{dudekV1996}).

\item[(ii)]

  Consider a $k$-ary group $(G,h)$ in which the $k$-ary operation $h$
  is distributive with respect to itself, i.e.,
 $$h(x_1^{i-1}, h(a_{1}^{k}), x_{i+1}^{k}) = h(h(x_1^{i-1},
  a_{1},x_{i+1}^{k}), \ldots , h(x_1^{i-1}, a_{k}, x_{i+1}^{k})),$$
 
 for all $1 \leq i \leq k$.  These type of groups are called
 \emph{autodistributive} $k$-ary groups (see
 Dudek~\cite{dudekauto1995}).

\end{itemize}

\end{remark}

\section{$(m,n)$-Semirings and Their Properties} \label{sem_sec}

\begin{definition} \label{def_gensemiring}

  An $(m,n)$-semiring is an algebraic structure
  $(\mathcal{R}, f, g)$ which satisfies the following axioms:

\begin{itemize}

\item[(i)] $(\mathcal{R},f)$ is an $m$-ary semigroup,

\item[(ii)] $(\mathcal{R},g)$ is an $n$-ary semigroup,

\item[(iii)] the $n$-ary operation $g$ is \emph{distributive} with
  respect to the $m$-ary operation $f$.

\end{itemize}

\end{definition}

\begin{example} \label{example_semiring}

  Let $\mathcal{B}$ be any Boolean algebra. Then $(\mathcal{B}, f,
  g)$ is an $(m, n)$-semiring where $f(A_{1}^{m})= A_1 \cup A_2\cup
  \ldots \cup A_m$ and $g(B_{1}^{n})=B_1 \cap B_2\cap \ldots \cap
  B_n$, for all $A_1,A_2,\ldots,A_m$ and $B_1,B_2,\ldots,B_n \in \mathcal{B}$.

\end{example}

In general, we have the following

\begin{theorem} \label{theo_init-sem}

  Let $(\mathcal{R}, +, \times)$ be an ordinary semiring.  Let $f$ be an
  $m$-ary operation and $g$ be an $n$-ary operation on $\mathcal{R}$
  as follows:
\[f(x_{1}^{m}) = \sum_{i=1}^m x_i, \qquad  \forall x_{1},x_{2},\ldots,x_{m}\in \mathcal{R},\]  
\[g(y_{1}^{n})=\prod_{i=1}^n y_i, \qquad  \forall y_{1},y_{2},\ldots,y_{n}\in \mathcal{R}.\]
Then $(\mathcal{R}, f, g)$ is an $(m, n)$-semiring. 

\end{theorem}

\begin{proof}

Omitted as obvious. \qedhere

\end{proof}

\begin{example} \label{example_sem1} 

  The following give us some $(m, n)$-semirings in different ways
  indicated by Theorem~\ref{theo_init-sem}.

\begin{itemize}

\item[(i)] Let $(\mathcal{R}, +, \times)$ be an ordinary semiring and
  $x_1,x_2,\ldots,x_n$ be in $\mathcal{R}$.  If we set
  \[g(x_1^n) = x_1 \times x_2 \times \ldots \times x_n,\]
   we get a $(2, n)$-semiring $(\mathcal{R}, +, g)$.

\item[(ii)] In an $(m, n)$-semiring $(\mathcal{R}, f, g)$, fixing
  elements $a_2^{m-1}$ and $b_2^{n-1}$, we obtain two binary
  operations as follows:

$$x \oplus y = f(x,a_2^{m-1}, y) \ \mathrm{and} \ x \otimes y = g(x, b_2^{n-1}, y).$$ 

Obviously, $(\mathcal{R}, \oplus, \otimes)$ is a semiring.

\item[(iii)] The set $\mathbb{Z}_{-}$ of all negative integers is not
  closed under the binary products, i.e., $\mathbb{Z}_{-}$ does not
  form a semiring, but it is a $(2, 3)$-semiring.

\end{itemize}

\end{example}

\begin{definition} \label{def_identityelem}

 Let $(\mathcal{R},f,g)$ be an $(m,n)$-semiring.  Then $m$-ary
 semigroup $(\mathcal{R},f)$ has an \emph{identity element}
 \(\mathbf{0}\) if
 
  \[x=f(\underbrace{\mathbf{0}, \ldots, \mathbf{0}}_{i-1}, x,
  \underbrace{\mathbf{0}, \ldots, \mathbf{0}}_{m-i})\] 

for all $x \in \mathcal{R}$ and $1 \leq i \leq m$.  We call
\(\mathbf{0}\) as an \emph{identity element} of $(m,n)$-semiring
$(\mathcal{R},f,g)$.

Similarly, $n$-ary semigroup $(\mathcal{R},g)$ has an \emph{identity
  element} \(\mathbf{1}\) if
 
  \[y=g(\underbrace{\mathbf{1}, \ldots,
    \mathbf{1}}_{j-1}, y, \underbrace{\mathbf{1}, \ldots,
    \mathbf{1}}_{n-j})\] for all $y \in \mathcal{R}$ and $1 \leq j
  \leq n$.

We call \(\mathbf{1}\) as an \emph{identity element} of
$(m,n)$-semiring $(\mathcal{R},f,g)$.

  We therefore call $\mathbf{0}$ the $f$-identity, and $\mathbf{1}$
  the $g$-identity.

\end{definition}

\begin{remark} \label{example_sem03}

  In an $(m,n)$-semiring $(\mathcal{R}, f, g)$, placing $\mathbf{0}$
  and $\mathbf{1}$, $(m-2)$ and $(n-2)$ times respectively, we obtain
  the following binary operations:

\begin{center}

  $x + y = f(x, \underbrace{\mathbf{0}, \ldots, \mathbf{0}}_{m-2}, y)$
  and $ x \times y =
  g(x,\underbrace{\mathbf{1},\ldots,\mathbf{1}}_{n-2},y), \qquad$ for
  all $x, y \in \mathcal{R}$.

\end{center}

\end{remark}

\begin{definition} 

  Let $(\mathcal{R}, f, g)$ be an $(m, n)$-semiring with an
  $f$-identity element $\mathbf{0}$ and $g$-identity element
  $\mathbf{1}$.  Then:

\begin{itemize}

\item[(i)]

  $\mathbf{0}$ is said to be \emph{multiplicatively absorbing} if it
  is absorbing in $(\mathcal{R}, g)$, i.e., if

  \[g(\mathbf{0}, x_{1}^{n-1}) = g(x_1^{n-1}, \mathbf{0}) =
  \mathbf{0}\] for all $x_1,x_2,\ldots,x_{n-1} \in \mathcal{R}$.

\item[(ii)] $(\mathcal{R}, f, g)$ is called
  \emph{zero-divisor free} if
  \[g(x_1, x_2, \ldots, x_n)=\mathbf{0}\] always implies $x_1 =
  \mathbf{0}$ or $x_2 = \mathbf{0}$ or $\ldots$ or $x_n = \mathbf{0}$.

  Elements $x_1, x_2, \ldots, x_{n-1} \in \mathcal{R}$ are called
  \emph{left zero-divisors} of $(m, n)$-semiring $(\mathcal{R}, f, g)$
  if there exists $a \neq \mathbf{0}$ and the following holds:

  \[g(x_1^{n-1}, a)=\mathbf{0}.\]

\item[(iii)] $(\mathcal{R}, f, g)$ is called
  \emph{zero-sum free} if
  \[f(x_1, x_2,\ldots, x_m)=\mathbf{0}\] always implies $x_1 = x_2 =
  \ldots = x_m=\mathbf{0}$.

\item[(iv)] $(\mathcal{R},f,g)$ is called \emph{additively
    cancellative} if the $m$-ary semigroup $(\mathcal{R},f)$ is
  cancellative, i.e.,

  \[f(x_1^{i-1},a,x_{i+1}^{m}) = f(x_1^{i-1},b,x_{i+1}^{m})
  \hspace{.2in} \Longrightarrow a=b\] for all $a,b,x_1,x_2,\ldots,x_{m} \in
  \mathcal{R}$ and for all $1\leq i \leq m$.

\item[(v)] $(\mathcal{R},f,g)$ is called \emph{multiplicatively
    cancellative} if the $n$-ary semigroup $(\mathcal{R},g)$ is
  cancellative, i.e.,
  \[g(x_1^{i-1},a,x_{i+1}^{n}) = g(x_1^{i-1},b,x_{i+1}^{n})
  \hspace{.2in} \Longrightarrow a=b\] for all $a,b,x_1,x_2,\ldots,x_{n} \in
  \mathcal{R}$ and for all $1\leq i \leq n$.

  Elements $x_1,x_2,\ldots,x_{n-1}$ are called \emph{left cancellable}
  in an $n$-ary semigroup \((\mathcal{R},g)\) if

  \[g(x_1^{n-1},a)=g(x_1^{n-1},b) \hspace{.2in} \Longrightarrow a=b\]
  for all $x_1,x_2,\ldots,x_{n-1},a,b \in \mathcal{R}.$

  $(\mathcal{R},f,g)$ is called \emph{multiplicatively left
    cancellative} if elements $x_1,x_2,\ldots,x_{n-1}$ $\in$ $\mathcal{R}
  \setminus$ $\{\mathbf{0}\}$ are multiplicatively left cancellable in
  $n$-ary semigroup $(\mathcal{R},g)$.

\end{itemize}  

\end{definition}

\begin{theorem} \label{left_cancellable}

  Let \((\mathcal{R},f,g)\) be an $(m,n)$-semiring with $f$-identity
  $\mathbf{0}$.

\begin{itemize}

\item[(i)] If elements $x_1,x_2,\ldots,x_{n-1} \in \mathcal{R}$ are
  multiplicatively left cancellable, then elements
  $x_1$,$x_2$,$\ldots$,$x_{n-1}$ are not left divisors.

\item[(ii)] If the $(m,n)$-semiring $(\mathcal{R},f,g)$ is
  multiplicatively left cancellative, then it is zero-divisor
  free.

\end{itemize} 

\end{theorem}


We have generalized Theorem~\ref{left_cancellable} from Theorem 4.4 of
Hebisch and Weinert~\cite{hebisch1998}.

We have generalized the definition of idempotents of
semirings given by Bourne~\cite{bourne1956} and Hebisch and
Weinert~\cite{hebisch1998}), as follows.

\begin{definition} \label{def_idem}
Let $(\mathcal{R},f,g)$ be an $(m,n)$-semiring.  Then: 

\begin{itemize}

\item[(i)] It is called \emph{additively idempotent} if
  $(\mathcal{R},f)$ is an idempotent $m$-ary semigroup, i.e., if
\[f(\underbrace{x,x,\ldots, x}_{m})=x\] for all $x \in \mathcal{R}$.

\item[(ii)] It is called \emph{multiplicatively idempotent} if
  $(\mathcal{R},g)$ is an idempotent $n$-ary semigroup, i.e., if

  \[g(\underbrace{y,y,\ldots ,y}_{n})=y\]
  for all $y \in \mathcal{R}$, $y \neq \mathbf{0}$.

\end{itemize}

\end{definition}

\begin{theorem} \label{thidempo}

  An $(m,n)$-semiring $(\mathcal{R},f,g )$ having at least two
  multiplicatively idempotent elements in the center is not multiplicatively
  cancellative.

\end{theorem}

\begin{proof} \label{}

  Let $a$ and $b$ be two multiplicatively idempotent elements in the center, $a \neq
  b$. Then:

  $g(\underbrace{\mathbf{1}, \ldots, \mathbf{1}}_{n-2}, a, b) =
  g(\underbrace{\mathbf{1}, \ldots, \mathbf{1}}_{n-2}, b, a)$

 which can be written as follows:

  $g(\underbrace{\mathbf{1}, \ldots, \mathbf{1}}_{n-2}, g(\stackrel{(n)}{a}), b) =
  g(\underbrace{\mathbf{1}, \ldots, \mathbf{1}}_{n-2}, g(\stackrel{(n)}{b}), a)$

  which is represented as:

  $g(\underbrace{\mathbf{1}, \ldots, \mathbf{1}}_{n-3}, g(\mathbf{1},\stackrel{(n-1)}{a}),a, b) =
  g(\underbrace{\mathbf{1}, \ldots, \mathbf{1}}_{n-3}, g(\mathbf{1},\stackrel{(n-1)}{b}), b, a)$.

  If the $(m,n)$-semiring $(\mathcal{R},f,g)$ is multiplicatively
  cancellative, then the following holds true:

$g(\underbrace{\mathbf{1}, \ldots, \mathbf{1}}_{n-3}, g(\mathbf{1},\stackrel{(n-1)}{a}), \mathbf{1}, \mathbf{1}) = g(\underbrace{\mathbf{1}, \ldots, \mathbf{1}}_{n-3}, g(\mathbf{1},\stackrel{(n-1)}{b}), \mathbf{1}, \mathbf{1})$,

$g(\mathbf{1},\stackrel{(n-1)}{a}) = g(\mathbf{1},\stackrel{(n-1)}{b})$,

which implies that $a = b$, which is a contradiction to the assumption
that $a \neq b$, therefore $(\mathcal{R},f,g)$ is not multiplicatively
cancellative. \qedhere

\end{proof}

We have generalized Exercise 2.7 in Chapter I of
Hebisch and Weinert~\cite{hebisch1998} to get the following.

\eat{
\begin{theorem} \label{add_idem_identity}

  An $(m,n)$-semiring $(R,f,g)$ having an $f$-identity element
  $\mathbf{0}$ is additively idempotent if and only if the following
  holds:

$$f(\underbrace{\mathbf{0},\ldots,\mathbf{0}}_{m})=\mathbf{0}.$$  

\end{theorem}
}

\begin{definition} \label{cong_def2} 

  Let ($\mathcal{R}, f ,g)$ be an $(m,n)$-semiring and $\sigma$ be an
  equivalence relation on $\mathcal{R}$. 
\begin{itemize}
\item[(i)]  
   Then $\sigma$ is called a
  \emph{congruence relation} or a \emph{congruence} of
  $(\mathcal{R},f,g)$, if it satisfies the following properties for all
  $1 \leq i \leq m$ and $1 \leq j \leq n$ :

\begin{itemize}

\item[(a)] if $x_i \sigma y_i$ then $f(x_1^m) \sigma f(y_1^m)$ and,

\item[(b)] if $z_j \sigma u_j$ then $g(z_1^n) \sigma g(u_1^n)$,

for all $x_1,x_2,\ldots,x_m$, $y_1,y_2,\ldots,y_m$ ,$z_1,z_2,\ldots,z_n$, $u_1,u_2,\ldots,u_n \in \mathcal{R}$.

\end{itemize}
\item[(ii)] Let $\sigma$ be a congruence on an algebra $\mathcal{R}$. Then the \emph{quotient} of $\mathcal{R}$ by $\sigma$, written as $\mathcal{R}/ \sigma$, is the algebra whose universe is $\mathcal{R}/ \sigma$ and whose fundamental operation satisfy 

$f^{\mathcal{R}/\sigma}(x_1,x_2,\ldots,x_m)$ = $ f^{\mathcal{R}}(x_1,x_2,\ldots,x_m)/ \sigma $

where $x_1,x_2,\ldots,x_m \in \mathcal{R}$~\cite{Burris1981}.
   
\end{itemize}  
  
\end{definition}



\begin{theorem} \label{theo_cong_rel}

  Let $(\mathcal{R}, f, g)$ be an $(m, n)$-semiring and the relation
  $\sigma$ be a congruence relation on $(\mathcal{R}, f, g)$. Then the
  quotient $(\mathcal{R} / \sigma, F, G)$ is an $(m,n)$-semiring under
  $F((x_1)/ \sigma, \ldots, (x_m)/\sigma)$ = $f(x_1^m )/\sigma$ and
  $G((y_1)/\sigma,\ldots, (y_n)/\sigma)$ =$ g(y_1^n)/\sigma$, for all
  $x_1,x_2,\ldots,x_m$ and $y_1,y_2,\ldots,y_n$ in $\mathcal{R}$.

\end{theorem}

\begin{proof}

Omitted as obvious. \qedhere

\end{proof}


\begin{definition} \label{def_homomorphism}

  We define homomorphism, isomorphism, and a product of two mappings
  as follows:

\begin{itemize}

\item[(i)]
  A mapping $\varphi: \mathcal{R} \rightarrow \mathcal{S}$ from
  $(m,n)$-semiring $(\mathcal{R},f,g)$ into $(m,n)$-semiring
  $(\mathcal{S},f',g')$ is called a \emph{homomorphism} if

\[\varphi(f(x_1^m))=f'(\varphi(x_1),\varphi(x_2),\ldots,\varphi(x_m))\] and 

\[\varphi(g(y_1^n))=g'(\varphi(y_1),\varphi(y_2),\ldots,\varphi(y_n))\]

for all $x_1,x_2,\ldots,x_m, y_1,y_2,\ldots,y_n \in \mathcal{R}$.

\item[(ii)] The $(m,n)$-semirings $(\mathcal{R},f,g)$ and
$(\mathcal{S},f',g')$ are called \emph{isomorphic} if there
exists one-to-one homomorphism from \(\mathcal{R}\) onto
\(\mathcal{S}\). One-to-one homomorphism is called \emph{isomorphism}.

\item[(iii)] If we apply mapping $\varphi:\mathcal{R} \rightarrow
  \mathcal{S}$ and then $\psi:\mathcal{S} \rightarrow \mathcal{T}$ on
  $x$ we get the mapping $(\psi \circ \varphi)(x)$ which is equal to
  $\psi(\varphi(x))$, where $x \in \mathcal{R}$. It is called the
  \emph{product} of $\psi$ and $\varphi$~\cite{hebisch1998}.

\end{itemize}

\end{definition}

We have generalized Definition~\ref{def_homomorphism} from Definition 2
of Allen~\cite{paulj1969}.

We have generalized the following theorem from Theorem 3.3 given by Hebisch and
Weinert~\cite{hebisch1998}.

\begin{theorem} \label{theo_homomorphism}

  Let $(\mathcal{R},f,g)$, $(\mathcal{S},f',g')$ and
  $(\mathcal{T},f'',g'')$ be $(m,n)$-semirings.  Then if 
  the following mappings $\varphi: (\mathcal{R},f,g) \rightarrow
  (\mathcal{S},f',g')$ and
 
  $\psi:(\mathcal{S},f',g') \rightarrow
  (\mathcal{T},f'',g'')$ are homomorphisms, then

  $\psi \circ \varphi:(\mathcal{R},f,g) \rightarrow
  (\mathcal{T},f'',g'')$ is also a homomorphism.

\end{theorem}

\begin{proof}
 Let $x_1,x_2,\ldots,x_m$ and $y_1,y_2,\ldots,y_n$ be in $\mathcal{R}$.  Then:
\begin{align*} 
 (\psi \circ \varphi)(f(x_1^m)) &= \psi(\varphi(f(x_1,x_2,\ldots,x_m))) \\ 
 &= \psi(f'(\varphi(x_1),\varphi(x_2),\ldots,\varphi(x_m))) \\
 &= f''(\psi(\varphi(x_1)),\psi(\varphi(x_2)),\ldots,\psi(\varphi(x_m))) \\
&= f''((\psi \circ \varphi)(x_1),(\psi \circ \varphi)(x_2),\ldots,(\psi \circ \varphi)(x_m)).
\end{align*}

In a similar manner, we can deduce that 
\[(\psi \circ \varphi)(g(y_1^n)) = g''((\psi \circ \varphi)(y_1),(\psi \circ \varphi)(y_2),\ldots,(\psi \circ \varphi)(y_n)). \]

Thus it is evident that $\psi \circ \varphi$ is a homomorphism from
\(\mathcal{R} \rightarrow \mathcal{T}\). \qedhere

\end{proof}

This proof is similar to that of Theorem 6.5 given by Burris and
Sankappanavar~\cite{Burris1981}.

\begin{definition} \label{def_kernel}

  Let $(\mathcal{R}, f, g)$ and $(\mathcal{S}, f', g')$ be
  $(m,n)$-semirings, and $\varphi: \mathcal{R} \rightarrow
  \mathcal{S}$ be a homomorphism.  Then the \emph{kernel} of
  $\varphi$, written as $\kr\varphi$ is defined as follows:

 \[\kr\varphi =\lbrace (a, b) \in \mathcal{R} \times \mathcal{R} \hspace{0.1in} \vert \hspace{0.1in} \varphi(a) =\varphi(b) \rbrace.\]
Generalization of Burris and Sankappanavar~\cite{Burris1981}.
\end{definition}

\begin{theorem} \label{cong_rel_theo1}

  Let $(\mathcal{R}, f, g)$ and $(\mathcal{S}, f', g')$ be
  $(m,n)$-semirings and $\varphi: \mathcal{R} \rightarrow
  \mathcal{S}$ be a homomorphism.  Then $\kr\varphi$ is a congruence
  relation on $\mathcal{R}$ and there exists a unique one-to-one
  homomorphism $\psi$ from $\mathcal{R}/ \kr \varphi$ into
  $\mathcal{S}$.

\end{theorem}

\begin{proof} 

Omitted as obvious. \qedhere

\end{proof}
\begin{corollary} \label{cong_col}

  Let $(\mathcal{R},f,g)$ be an $(m,n)$-semiring and $\rho$ and
  $\sigma$ be congruence relations on $\mathcal{R}$, with $\rho
  \subseteq \sigma$. Then $\sigma/ \rho= \lbrace \rho(x),\rho(y)
  \hspace{0.1in} \vert \hspace{0.1in} (x,y) \in \sigma \rbrace$ is a
  congruence relation on $\mathcal{R}/ \rho$, and $(\mathcal{R}/
  \rho)/(\sigma/ \rho) \cong \mathcal{R}/ \sigma$.

\end{corollary}

 \begin{lemma} \label{zerogen}

   Let $x_1, x_2, \ldots, x_m, y_1, y_2, \ldots, y_n \in
   \mathcal{R}$. Then:

\begin{itemize}

\item[(i)]
$\underbrace{f(f(\ldots f(f}_{m}(x_1, \underbrace{\mathbf{0},\ldots,\mathbf{0}}_{m-1}),x_2,\underbrace{\mathbf{0},\ldots,\mathbf{0}}_{m-2}),\ldots),x_m,\underbrace{\mathbf{0},\ldots,\mathbf{0}}_{m-2})$ = $f(x_1,x_2,\ldots,x_m)$,

\item[(ii)]
$\underbrace{g(g(\ldots g(g}_{n}(y_1,\underbrace{\mathbf{1},\ldots,\mathbf{1}}_{n-1}),y_2,\underbrace{\mathbf{1},\ldots,\mathbf{1}}_{n-2}),\ldots),y_n,\underbrace{\mathbf{1},\ldots,\mathbf{1}}_{n-2})$ = $g(y_1,y_2,\ldots,y_n).$

\end{itemize}
\end{lemma}

\begin{proof}

\begin{itemize}

\item[(i)]
\begin{equation} \label{zerogeneqn1} 
\underbrace{f(f(\ldots f(f}_{m}(x_1,\underbrace{\mathbf{0},\ldots,\mathbf{0}}_{m-1}),x_2,\underbrace{\mathbf{0},\ldots,\mathbf{0}}_{m-2}),\ldots),x_m,
\underbrace{\mathbf{0},\ldots,\mathbf{0}}_{m-2}). 
\end{equation}

By associativity (Definition~\ref{definition_additive} (i)),
\eqref{zerogeneqn1} is equal to

\[\underbrace{f(f(\ldots f(f}_{m}(\mathbf{0},\underbrace{\mathbf{0},\ldots,\mathbf{0}}_{m-1}),x_1,x_2,\underbrace{\mathbf{0},\ldots,\mathbf{0}}_{m-3}),\ldots),x_m,\underbrace{\mathbf{0},\ldots,\mathbf{0}}_{m-2})\] 
\[ = \underbrace{f(f(\ldots f(f}_{m-1}(\mathbf{0},x_1,x_2,\underbrace{\mathbf{0},\ldots,\mathbf{0}}_{m-3}),x_3,\underbrace{\mathbf{0},\ldots,\mathbf{0}}_{m-2}),\ldots),x_m,\underbrace{\mathbf{0},\ldots,\mathbf{0}}_{m-2}) \]
\[= \underbrace{f(f(\ldots f(f}_{m-1}(\underbrace{\mathbf{0},\ldots,\mathbf{0}}_{m}),x_1,x_2,x_3,\underbrace{\mathbf{0},\ldots,\mathbf{0}}_{m-4}),\ldots),x_m,\underbrace{\mathbf{0},\ldots,\mathbf{0}}_{m-2})\]
\[ = \underbrace{f(f(\ldots f(f}_{m-2}(\mathbf{0},x_1,x_2,x_3,\underbrace{\mathbf{0},\ldots,\mathbf{0}}_{m-4}),\ldots),x_m,\underbrace{\mathbf{0},\ldots,\mathbf{0}}_{m-2})\]
\hspace{2in} \vdots
\[= f(f(x_1,x_2,\ldots,x_{m-1},\mathbf{0}),x_m,\underbrace{\mathbf{0},\ldots,\mathbf{0}}_{m-2})\] 
\[= f(f(x_1,x_2,\ldots,x_{m-1},x_m),\underbrace{\mathbf{0},\ldots,\mathbf{0}}_{m-1})\] 
\[= f(x_1,x_2,\ldots,x_m).\]
\item[(ii)]
Similar to part (i). \qedhere
\end{itemize} 
\end{proof}

\section{ Partial Ordering On Fault Tolerance} \label{partial_ordering}

In this Section we use $x_i,y_i$, etc., where $i \in \mathbb{Z}_+$ to
denote individual system components that are assumed to be
\emph{atomic} at the level of discussion, i.e., they have no
components or sub-systems of their own.  We use \emph{component} to
refer to such an atomic part of a system, and \emph{subsystem} to
refer to a part of a system that is not necessarily atomic.  We assume
that components and subsystems are disjoint, in the sense that if
fail, they fail independently and do not affect the functioning of
other components.

Let $\mathcal{U}$ be a universal set of all systems in the domain of
discourse as given by Rao~\cite{SRaoAlgeb2008}, and let $f$ be a
mapping $f:\mathcal{U}^m \to \mathcal{U}$, i.e., $f$ is an $m$-ary
operation.  Likewise, let $g$ be an $n$-ary operation.

\begin{definition} \label{f_def}

We define $f$ and $g$ operations for systems as follows:

\begin{itemize} 
\item[(i)]
  $f$ is an $m$-ary operation which applies on systems made up of $m$
  components or subsystems, where if any one of the
  components or subsystems fails, then the whole system fails.

  Let a system made up of $m$ components $x_1,x_2, \ldots, x_m$, then
  the system over operation $f$ is represented as $f(x_1,x_2, \ldots,
  x_m)$ for all $x_1, x_2, \ldots, x_m \in \mathcal{U}$. The system
  $f(x_1, x_2, \ldots, x_m)$ fails when any of the components $x_1,
  x_2, \ldots, x_m$ fails.

\item[(ii)] $g$ is an $n$-ary operation which applies on a system
  consisting of $n$ components or subsystems, which fails if all the
  components or subsystems fail; otherwise it continues working even if
  a single component or subsystem is working properly.

  Let a system consist of $n$ components $x_1,x_2 , \ldots , x_n$,
  then the system over operation $g$ is represented as $g(x_1,x_2 ,
  \ldots , x_n)$ for all $x_1,x_2 , \ldots , x_n \in \mathcal{U}$.
  The system $g(x_1,x_2 , \ldots , x_n)$ fails when all the components
  $x_1,x_2 , \ldots , x_n$ fail.

\end{itemize}
\end{definition} 

Consider a partial ordering relation \(\preccurlyeq\) on
\(\mathcal{U}\), such that \((\mathcal{U}, \preccurlyeq)\) is a
partially ordered set (poset).  This is a \emph{fault-tolerance
  partial ordering} where $f(x_{1}^{m}) \preccurlyeq f(y_{1}^{m})$
means that $f(x_{1}^{m})$ has a lower measure of some fault metric
than $f$($y_{1}^{m}$) and $f$($x_{1}^{m}$) has a better fault
tolerance than $f(y_{1}^{m})$, for all $f(x_{1}^{m}), f(y_{1}^{m}) \in
\mathcal{U}$ (see Rao~\cite{SRao2008} for more details) and
$x_{1},x_{2},\ldots,x_{m}, y_{1},y_{2},\ldots,y_{m}$ are disjoint components.

Assume that \(\mathbf{0}\) represents the atomic system ``which is
always up'' and \(\mathbf{1}\) represents the system ``which is always
down'' (see Rao~\cite{SRao2008}).


\begin{observation} \label{obs_identity}
We observe the following for all disjoint components
 $x_1$,$x_2,\ldots,x_m$, $y_1,y_2,\dots,y_n$, which are in $\mathcal{U}$:

\begin{itemize}

\item[(i)] \(g(y_1^{j-1}, \mathbf{0}, y_{j+1}^{n}) = \mathbf{0}\) for
  all $1\leq j \leq n$.

  This is so since \(\mathbf{0}\) represents the component or system
  which never fails, and as per the definition of $g$, the system as a
  whole fails if all the components fail, and otherwise it continues
  working even if a single component is working properly.  In a system
  $g(y_1^{j-1}, \mathbf{0}, y_{j+1}^{n})$, even if all other
  components $y_1^{j-1}$ and $y_{j+1}^{n}$ fail even then
  \(\mathbf{0}\) is up and the system is always up.

\item[(ii)] \(f(x_1^{i-1}, \mathbf{1}, x_{i+1}^{m})= \mathbf{1}\) for
  all $1\leq i \leq m$.

  This is so since \(\mathbf{1}\) represents the component or system
  which is always down, and as per the definition of $f$ if either of
  the component fails, then the whole system fails. Thus, even though
  all other components are working properly but due to the component
  \(\mathbf{1}\) the system is always down.

\end{itemize}

\end{observation}

\begin{definition} \label{defpartialordering}

  If $(\mathcal{U},f,g)$ is an $(m,n)$-semiring and \((\mathcal{U},
  \preccurlyeq)\) is a poset, then \((\mathcal{U}, f, g,
  \preccurlyeq)\) is a \emph{partially ordered $(m,n)$-semiring} if
  the following conditions are satisfied for all $x_{1},x_{2},\ldots,x_{m},
  y_{1},y_{2},\ldots,y_{n}, a, b \in \mathcal{U}$ and $1 \leq i \leq m , 1 \leq j
  \leq n$.

\begin{itemize}

\item[(i)] If $a \preccurlyeq b$, then $f(x_1^{i-1}, a, x_{i+1}^{m})
  \preccurlyeq f(x_1^{i-1}, b, x_{i+1}^{m})$.

\item[(ii)] If $a \preccurlyeq \textit{b}$, then $g(y_1^{j-1}, a,
  y_{j+1}^{n}) \preccurlyeq g(y_1^{j-1}, b, y_{j+1}^{n})$.

\end{itemize}

\end{definition}

\begin{remark} \label{remarks_1}

  As it is assumed that $\mathbf{0}$ is the system which is always up,
  it is more fault tolerant than any of the other systems or
  components. Therefore $\mathbf{0} \preccurlyeq a$, for all $a \in
  \mathcal{U}$. Similarly, $a \preccurlyeq \mathbf{1}$ because $\mathbf{1}$ is the system that always fails and therefore it is the least fault tolerant; every other system
  is more fault-tolerant than it.

\end{remark}

\begin{observation} \label{observ1}

The following are obtained for all disjoint components $r,s,x_i,y_j$,
$a_i,b_j$, which are in $\mathcal{U}$, where $1 \leq i \leq m$ , $1 \leq j \leq n$:

\begin{itemize}

\item[(i)]

  $\mathbf{0} \preccurlyeq f(x_1^{i-1} , r , x_{i+1}^{m}) \preccurlyeq
  \mathbf{1}$.

\item[(ii)]

  $\mathbf{0} \preccurlyeq g(y_1^{j-1} , s , y_{j+1}^{n}) \preccurlyeq
  \mathbf{1}$.

\item[(iii)]

  $\mathbf{0} \preccurlyeq g(y_1^{j-1} , f(a_{1}^m) , y_{j+1}^{n})
  \preccurlyeq \mathbf{1}$.

\item[(iv)]

  $\mathbf{0} \preccurlyeq f(x_1^{i-1} , g(b_{1}^n) , x_{i+1}^{m})
  \preccurlyeq \mathbf{1}$.

\end{itemize}

\end{observation}

From the above description of $\mathbf{0}$ and $\mathbf{1}$, the
observation is quite obvious. Case (i) shows that \(\mathbf{0}\) is
less faulty than $f(x_1^{i-1}, r, x_{i+1}^{m})$, and $f(x_1^{i-1}, r,
x_{i+1}^{m})$ is less faulty than \(\mathbf{1}\).  Similarly, case
(ii) shows that \(\mathbf{0}\) is more fault-tolerant than
$g(y_1^{j-1}, s, y_{j+1}^{n})$ and $g(y_1^{j-1}, s, y_{j+1}^{n})$ is
more fault-tolerant than \(\mathbf{1}\). Likewise, case (iii) shows
the operation $g$ over $y_1^{j-1}$, $y_{j+1}^{n}$ and $f$ of $a_1^m$
to be less faulty than \(\mathbf{1}\) and more faulty than
\(\mathbf{0}\), and a similar interpretation is made for (iv).

\begin{lemma} \label{lemma2}

  If $\preccurlyeq$ is a fault-tolerance partial order and $x_1,x_2,\ldots,x_m$,
  $y_1,y_2,\ldots,y_m$, $z_1,z_2,\ldots,z_n$,$u_1,u_2,\ldots,u_n$ are 
  disjoint components, which are in $\mathcal{U}$, where $m,n \in \mathbb{Z_+}$, then for all $1 \leq i \leq m$ and $1 \leq j \leq n$ following holds true:

\begin{itemize}

\item[(i)] if $x_i \preccurlyeq y_i$, then $f(x_{1}^m) \preccurlyeq
  f(y_{1}^m)$, and, 

\item[(ii)] if $z_j \preccurlyeq u_j$, then $g(z_{1}^n) \preccurlyeq
  g(u_{1}^n)$.

\end{itemize}

\end{lemma}

\begin{proof}

\begin{itemize}

\item[(i)] Since $x_i \preccurlyeq y_i$ for all $1 \leq i \leq m$,
  we have:

\[x_1 \preccurlyeq y_1\]
which is represented as follows:
\begin{equation} \label{pleqq1}
 f(\underbrace{\mathbf{0},\ldots,\mathbf{0}}_{m-1},x_1) \preccurlyeq f(\underbrace{\mathbf{0},\ldots,\mathbf{0}}_{m-1},y_1)
\end{equation}
and

\begin{equation} \label{pleqq2}
f(\underbrace{\mathbf{0},\ldots,\mathbf{0}}_{m-1},x_2) \preccurlyeq f(\underbrace{\mathbf{0},\ldots,\mathbf{0}}_{m-1},y_2).
\end{equation}

By $f$ operation on both sides of \eqref{pleqq1} with $y_2$, we get:

\begin{equation} \label{pleqq3}
f(f(\underbrace{\mathbf{0} , \ldots ,  \mathbf{0}}_{m-1}, x_1) , y_2,\underbrace{\mathbf{0} , \ldots ,  \mathbf{0}}_{m-2}) \preccurlyeq f(f(\underbrace{\mathbf{0} , \ldots ,  \mathbf{0}}_{m-1}, y_1), y_2,\underbrace{\mathbf{0} , \ldots ,  \mathbf{0}}_{m-2}).
\end{equation}

By $f$ operation on both sides of \eqref{pleqq2} with $x_1$:

\begin{equation} \label{pleqq4}
f(f(\underbrace{\mathbf{0} , \ldots ,  \mathbf{0}}_{m-1}, x_2),x_1,\underbrace{\mathbf{0} , \ldots ,  \mathbf{0}}_{m-2} ) \preccurlyeq f(f(\underbrace{\mathbf{0} , \ldots ,  \mathbf{0}}_{m-1}, y_2) , x_1 ,\underbrace{\mathbf{0} , \ldots ,  \mathbf{0}}_{m-2}).
\end{equation}

From \eqref{pleqq3} and \eqref{pleqq4}, we get:
\[f(f(\underbrace{\mathbf{0} , \ldots ,  \mathbf{0}}_{m-1}, x_1) , y_2,\underbrace{\mathbf{0} , \ldots ,  \mathbf{0}}_{m-2} ) \preccurlyeq f(f(\underbrace{\mathbf{0} , \ldots ,  \mathbf{0}}_{m-1}, y_1) , y_2,\underbrace{\mathbf{0} , \ldots ,  \mathbf{0}}_{m-2}).\]

Similarly, we find for $m$ terms:

\begin{equation} \label{pleqq6}
\begin{split}
  \underbrace{f(\ldots(f(f}_{m}(\underbrace{\mathbf{0} , \ldots ,  \mathbf{0}}_{m-1},x_1) , x_2,\underbrace{\mathbf{0} , \ldots ,  \mathbf{0}}_{m-2}) , \ldots ), x_m,\underbrace{\mathbf{0} , \ldots ,  \mathbf{0}}_{m-2})
\\
\preccurlyeq \underbrace{f(\ldots(f(f}_{m}(\underbrace{\mathbf{0} , \ldots ,  \mathbf{0}}_{m-1},y_1) , y_2,\underbrace{\mathbf{0} , \ldots ,  \mathbf{0}}_{m-2}) , \ldots ), y_m,\underbrace{\mathbf{0} , \ldots ,  \mathbf{0}}_{m-2}).
\end{split}
\end{equation}

From Lemma~\ref{zerogen}, \eqref{pleqq6} may be represented as

\[ f(x_{1},x_{2},\ldots,x_{m}) \preccurlyeq f(y_{1},y_{2},\ldots,y_{m})\]
so
 \[ f(x_{1}^m) \preccurlyeq f(y_{1}^m). \]


\item[(ii)]
Since $z_j \preccurlyeq y_j$, for all $1 \leq j \leq n$
\[g(\underbrace{\mathbf{1},\ldots,\mathbf{1}}_{n-1},z_1) \preccurlyeq g(\underbrace{\mathbf{1},\ldots,\mathbf{1}}_{n-1},u_1)\]
and 

\[g(\underbrace{\mathbf{1},\ldots,\mathbf{1}}_{n-1},z_2) \preccurlyeq g(\underbrace{\mathbf{1},\ldots,\mathbf{1}}_{n-1},u_2). \]


\eat{

by $g$ operation on both sides of \eqref{lemma20} with $y_2$, we get the following

\begin{equation} \label{lemma22}
g(g(\underbrace{\mathbf{1} , \ldots ,  \mathbf{1}}_{n-1}, x_1) , y_2,\underbrace{\mathbf{1} , \ldots ,  \mathbf{1}}_{n-2}) \preccurlyeq g(g(\underbrace{\mathbf{1} , \ldots ,  \mathbf{1}}_{n-1}, y_1), y_2,\underbrace{\mathbf{1} , \ldots ,  \mathbf{1}}_{n-2})
\end{equation}

by $g$ operation on both sides of \eqref{lemma21} with $x_1$, we get

\begin{equation} \label{lemma23}

g(g(\underbrace{\mathbf{1} , \ldots ,  \mathbf{1}}_{n-1}, x_1) , x_2,\underbrace{\mathbf{1} , \ldots ,  \mathbf{1}}_{n-2}) \preccurlyeq g(g(\underbrace{\mathbf{1} , \ldots ,  \mathbf{1}}_{n-1}, y_2), x_1,\underbrace{\mathbf{1} , \ldots ,  \mathbf{1}}_{n-2}).

\end{equation}

From \eqref{lemma22} and \eqref{lemma23}, we get the following 

\[g(g(\underbrace{\mathbf{1} , \ldots ,  \mathbf{1}}_{n-1}, x_1) , x_2,\underbrace{\mathbf{1} , \ldots ,  \mathbf{1}}_{n-2} ) \preccurlyeq g(g(\underbrace{\mathbf{1} , \ldots ,  \mathbf{1}}_{n-1}, y_1) , y_2,\underbrace{\mathbf{1} , \ldots ,  \mathbf{1}}_{n-2}).\]
}


After following similar steps as seen in part (i), we use the $g$
operation for $n$ terms,

\[\begin{split} 
\underbrace{g(\ldots(g(g}_{n}(\underbrace{\mathbf{1} , \ldots ,  \mathbf{1}}_{n-1},z_1) , z_2,\underbrace{\mathbf{1} , \ldots ,  \mathbf{1}}_{n-2}) , \ldots ), z_n,\underbrace{\mathbf{1} , \ldots ,  \mathbf{1}}_{n-2}) 
\\
\preccurlyeq \underbrace{g(\ldots(g(g}_{n}(\underbrace{\mathbf{1} , \ldots ,  \mathbf{1}}_{n-1},u_1) , u_2,\underbrace{\mathbf{1} , \ldots ,  \mathbf{1}}_{n-2}) , \ldots ), u_n,\underbrace{\mathbf{1} , \ldots ,  \mathbf{1}}_{n-2}) \end{split}\]
which is represented as
\[ g(z_1 , z_2 , \ldots , z_n) \preccurlyeq g(u_1 , u_2 , \ldots , u_n)\]
and so
\[ g(z_{1}^n) \preccurlyeq g(u_{1}^n).\qedhere   \] 

\end{itemize}

\end{proof}

\begin{theorem} \label{posetth}

  If $\preccurlyeq$ is a fault-tolerance partial order and given
  disjoint components $a_i, c_j, b_i, d_j$ in $\mathcal{U}$, where $1
  \leq i \leq m$, $1 \leq j \leq n$ and $1 \leq k \leq m$, the following obtain:

\begin{itemize}

\item[(i)] If $a_{i} \preccurlyeq b_{i}$, then:

  $g(y_1^{j-1} , f(a_{1}^m) , y_{j+1}^{n}) \preccurlyeq g(y_1^{j-1} ,
  f(b_{1}^m),y_{j+1}^{n})$, for all $y_1,y_2,\ldots,y_n \in \mathcal{U}$.

\item[(ii)] If $c_{j} \preccurlyeq d_{j}$, then:

  $f(x_1^{k-1} , g(c_{1}^n) , x_{k+1}^{m}) \preccurlyeq f(x_1^{k-1} ,
  g(d_{1}^n),x_{k+1}^{m})$, for all $x_1,x_2,\ldots,x_m \in \mathcal{U}$.

\end{itemize}

\end{theorem}

\begin{proof}

\begin{itemize}

\item[(i)] Since $a_{i} \preccurlyeq b_{i}$, $\qquad$ for all $1 \leq
  i \leq m$.

Therefore, from Lemma~\ref{lemma2} (i) 
\[f(a_{1}^m) \preccurlyeq f(b_{1}^m), \qquad \forall a_1,a_2,\ldots,a_m,b_1,b_2,\ldots,b_m \in \mathcal{U}.\]

From Definition~\ref{defpartialordering} of a partially ordered
$(m,n)$-semiring, we deduce that
\[g(y_{1}^{j-1}, f(a_{1}^m) , y_{j+1}^{n}) \preccurlyeq g(y_{1}^{j-1} , f(b_{1}^m) , y_{j+1}^{n})\]
for all $1 \leq j \leq n$.

\item[(ii)] Since $c_{j} \preccurlyeq d_{j}$, $\qquad$ for all $1 \leq
  j \leq n$, from Lemma \ref{lemma2} (ii), we find that

\[g(c_{1}^n) \preccurlyeq g(d_{1}^n),  \qquad \forall c_1,c_2,\ldots,c_n,d_1,d_2,\ldots,d_n \in \mathcal{U}.\]

From Definition \ref{defpartialordering} of a partially ordered
$(m,n)$-semiring, we deduce that

\[f(x_{1}^{k-1},g(c_{1}^n),x_{k+1}^{m}) \preccurlyeq f(x_{1}^{k-1},g(d_{1}^n),x_{k+1}^{m})\]
for all $1 \leq k \leq m$. \qedhere

\end{itemize}

\end{proof}

\begin{lemma} \label{lemmax1}

  If $\preccurlyeq$ is a fault-tolerance partial order and $x_i , y_j$
  are disjoint components which are in $\mathcal{U}$, where $1 \leq i
  \leq m$ and $1 \leq j \leq n$, we get the following:

\begin{itemize}

\item[(i)] $x_i \preccurlyeq f(x_{1},x_{2},\ldots,x_{m})$,  

\item[(ii)] $g(y_{1},y_{2},\ldots,y_{n}) \preccurlyeq y_j$.

\end{itemize}

\end{lemma}

\begin{proof}

\begin{itemize}
\item[(i)] As
\begin{equation} \label{lemxeqn1}
 \mathbf{0} \preccurlyeq x_{1},
\end{equation} 
by $f$ operation on both sides of \eqref{lemxeqn1} with  $x_{i}$, we get
\[f(\mathbf{0},x_{i},\underbrace{\mathbf{0} , \ldots , \mathbf{0}}_{m-2}) \preccurlyeq f(x_{1},x_{i},\underbrace{\mathbf{0} , \ldots , \mathbf{0}}_{m-2}).\]
Therefore,
\[ x_{i} \preccurlyeq f(x_{1},x_{i},\underbrace{\mathbf{0} , \ldots , \mathbf{0}}_{m-2}).\]

Similarly, we obtain:

\begin{equation} \label{lemxeqn4}
\begin{split}
 x_{i} \preccurlyeq f(x_{1},x_{i},\underbrace{\mathbf{0} , \ldots , \mathbf{0}}_{m-2})
\preccurlyeq \ldots \preccurlyeq f(x_{1}, x_{2},x_{i}, \ldots, x_{m-1},\mathbf{0}) 
\\
\preccurlyeq f(x_{1},x_{2},\ldots,x_{m}).
\end{split}
\end{equation} 

Hence,
\[ x_{i} \preccurlyeq f(x_{1},x_{2},\ldots,x_{m})\]
for all $1 \leq i \leq m$.

\item[(ii)] As

\begin{equation} \label{lemxeqn11}
  y_{1}  \preccurlyeq \mathbf{1},
\end{equation} 

by $g$ operation on both sides of \eqref{lemxeqn11} with  $y_{j}$, we get 
\[ g(y_{1},y_{j},\underbrace{\mathbf{1} , \ldots , \mathbf{1}}_{n-2}) \preccurlyeq y_j.\]

Similarly, we obtain:

\begin{equation} \label{lemxeqn14}
\begin{split}
g(y_{1},y_{2},\ldots,y_{n}) \preccurlyeq g(y_{1}, y_{2},y_{j}, \ldots, y_{n-1},\mathbf{1})  \preccurlyeq  \ldots \preccurlyeq 
\\
g(y_{1},y_{j},\underbrace{\mathbf{1} , \ldots , \mathbf{1}}_{n-2}) \preccurlyeq y_{j}.
\end{split}
\end{equation} 

Hence,

  \[g(y_{1},y_{2},\ldots,y_{n}) \preccurlyeq y_{j}\] 
for all $1 \leq j \leq n$. \qedhere
\end{itemize}

\end{proof}

\begin{corollary} \label{corollaryx2}

  If $\preccurlyeq$ is a fault-tolerance partial order, then the
  following hold for all disjoint components $x_i ,y_j$ which are
  elements of $\mathcal{U}$, where $1 \leq i \leq m$, $1 \leq j \leq
  n$ and $k, t \in \mathbb{Z_+}$:

\begin{itemize}

\item[(i)] $f(x_1, x_2, \ldots, x_k,\underbrace{\mathbf{0} , \ldots ,
    \mathbf{0}}_{m-k}) \preccurlyeq f(x_{1}^m)$,

  where $k < m$; and

\item[(ii)] $g(y_{1}^n) \preccurlyeq g(y_1, y_2, \ldots,
  y_t,\underbrace{\mathbf{1} , \ldots , \mathbf{1}}_{n-t})$,

  where $t < n$.

\end{itemize}

\end{corollary}
\begin{proof}
\begin{itemize}

\item[(i)] From \eqref{lemxeqn4} we deduce that, 

\[\begin{split}
f(x_{1},\ldots ,x_{k},\underbrace{\mathbf{0} , \ldots , \mathbf{0}}_{m-k})
\preccurlyeq f(x_{1}, \ldots, x_{k+1},\underbrace{\mathbf{0} , \ldots , \mathbf{0}}_{m-k-1})
\\
 \preccurlyeq \ldots \preccurlyeq f(x_{1},x_{2},\ldots,x_{m}).
\end{split}\]

Therefore,
\[ f(x_{1},\ldots ,x_{k},\underbrace{\mathbf{0} , \ldots , \mathbf{0}}_{m-k})
\preccurlyeq f(x_{1}^{m}). \]

\item[(ii)] As in part (i), we deduce from \eqref{lemxeqn14} that:

\[g(y_{1}^n)  \preccurlyeq g(y_1, y_2, \ldots, y_t,\underbrace{\mathbf{1} , \ldots , \mathbf{1}}_{n-t}).\qedhere \] 

\end{itemize}
\end{proof}

\eat{
\begin{theorem} \label{congnewth1}

  If $\preccurlyeq$ is a fault-tolerance partial order and $x_1,x_2,\ldots,x_m$,$y_1,y_2,\ldots,y_m$ ,
  $z_1,z_2,\ldots,z_n$, $u_1,u_2,\ldots,u_n$ are disjoint components in $\mathcal{U}$, then
  following hold:

\begin{itemize}

\item[(i)] If $f(x_i, \ldots, x_m,\underbrace{\mathbf{0} , \ldots ,
    \mathbf{0}}_{i-1}) \preccurlyeq f(y_i, \ldots,
  y_m,\underbrace{\mathbf{0} , \ldots , \mathbf{0}}_{i-1})$,

  then $f(x_1,x_2,\ldots,x_m) \preccurlyeq f(y_1,y_2,\ldots,y_m)$ for
  all $1 < i < m$.

\item[(ii)] If $g(z_j, \ldots, z_n,\underbrace{\mathbf{1} , \ldots ,
    \mathbf{1}}_{j-1}) \preccurlyeq g(u_j, \ldots,
  u_n,\underbrace{\mathbf{1} , \ldots , \mathbf{1}}_{j-1})$,

  then $g(z_1,z_2,\ldots,z_n) \preccurlyeq g(u_1,u_2,\ldots,u_n)$ for
  all $1 < j < n$.

\end{itemize}

\end{theorem}
}

$f(\stackrel{(m)}{f(a_{1}^m)})$ represents the system which is
obtained after applying the $f$ operation on $m$ repeated $f(a_{1}^m)$
systems or subsystems.  Similarly, $g(\stackrel{(n)}{g(b_{1}^n)})$
represents the system which is obtained after applying the $g$
operation on $n$ repeated $g(b_{1}^n)$ systems or subsystems.

\begin{theorem} \label{theorem40}

  If $\preccurlyeq$ is a fault-tolerance partial order, and components
  $x_{1},x_{2},\ldots,x_m$ , $ y_{1},y_{2},\ldots,y_n$ are disjoint components and are in $\mathcal{U}$, then:

\begin{itemize}

\item[(i)] $f(x_{1}^m) \preccurlyeq f(\stackrel{(m)}{f(x_{1}^m)})$ ,

\item[(ii)] $g(\stackrel{(n)}{g(y_{1}^n)}) \preccurlyeq g(y_{1}^n)$.

\end{itemize}

\end{theorem}

\begin{corollary} \label{corollary3} 

  The following hold for all disjoint components $x_1,\ldots,x_m$,
  
  $ z_1,\ldots,z_n$, $ y_1,\ldots,y_m$, $u_1,\ldots,u_n$, which are elements of $\mathcal{U}$, where $m,n \in \mathbb{Z_{+}}$:

\begin{itemize}
\item[(i)] If $f(\stackrel{(m)}{f(x_{1}^m)}) \preccurlyeq f(y_{1}^m)$, then 

$f(x_{1}^m) \preccurlyeq f(y_{1}^m)$.

\item[(ii)] If $g(z_{1}^n) \preccurlyeq g(\stackrel{(n)}{g(u_{1}^n)})$, then

$g(z_{1}^n) \preccurlyeq g(u_{1}^n)$.
\end{itemize}

\end{corollary}

\begin{proof}

\begin{itemize}

\item[(i)] $f(\stackrel{(m)}{f(x_{1}^m)}) \preccurlyeq f(y_{1}^m)$

  and from Theorem~\ref{theorem40}, $f(x_{1}^m) \preccurlyeq
  f(\stackrel{(m)}{f(x_{1}^m)})$.

  Therefore, $f(x_{1}^m) \preccurlyeq f(y_{1}^m)$.

\item[(ii)] The proof is very similar to that of part (i). \qedhere

\end{itemize}
\end{proof}

\begin{corollary} \label{corollary4}

  Let $k$ and $t$ be positive integers and $k < m$, $t < n$. Given
  disjoint components $x_1,\ldots,x_m$, $y_1,y_2,\ldots,y_m$, 
  
  $z_1,z_2,\dots,z_n$,$u_1,u_2,\ldots,u_n$ that are in $\mathcal{U}$, the following hold:

\begin{itemize}

\item[(i)] If
  $f(\underbrace{\mathbf{0},\ldots,\mathbf{0}}_{m-k},\stackrel{(k)}{f(x_{1}^m)})
  \preccurlyeq f(y_{1}^m)$, then $f(x_{1}^m)
  \preccurlyeq f(y_{1}^m)$.

\item[(ii)] If $g(z_{1}^n) \preccurlyeq
  g(\underbrace{\mathbf{1},\ldots,\mathbf{1}}_{n-t},\stackrel{(t)}{g(u_{1}^n)})$, then $g(z_{1}^n) \preccurlyeq g(u_{1}^n)$.

\end{itemize}

\end{corollary}

\begin{proof}
Similar to Corollary~\ref{corollary3}. \qedhere
\end{proof}

\begin{theorem} \label{thgen}

  Let $\preccurlyeq$ be a fault-tolerance partial order and $x_i \preccurlyeq y_i$ and $z_j \preccurlyeq u_j$ for all $x_{i}, y_{i},z_{j}, u_{j} \in \mathcal{U}$, where $1 \leq i \leq m$ and $1 \leq j \leq n$.  Then the following obtain:

\begin{itemize}

\item[(i)]
  $f(\stackrel{(m)}{f(x_{1}^m)}) \preccurlyeq
  f(\stackrel{(m)}{f(y_{1}^m)})$,

\item[(ii)]
  $g(\stackrel{(n)}{g(z_{1}^n)}) \preccurlyeq
  g(\stackrel{(n)}{g(u_{1}^n)})$,

\item[(iii)] 
  $f(\stackrel{(m)}{g(z_{1}^n)}) \preccurlyeq
  f(\stackrel{(m)}{g(u_{1}^n)})$,

\item[(iv)] 
  $g(\stackrel{(n)}{f(x_{1}^m)}) \preccurlyeq
  g(\stackrel{(n)}{f(y_{1}^m)})$.

\end{itemize}
\end{theorem}

\begin{proof}

\begin{itemize}

\item[(i)] As
\[x_i \preccurlyeq y_i, \qquad 1 \leq i \leq m,\]
from Lemma \ref{lemma2} (i), we get
\[f(x_1^m) \preccurlyeq f(y_1^m).\]

This is written as
\begin{equation} \label{thgeneqn110}
f(\underbrace{\mathbf{0} , \ldots ,  \mathbf{0}}_{m-1}, f(x_{1}^m)) \preccurlyeq f(\underbrace{\mathbf{0} , \ldots ,  \mathbf{0}}_{m-1}, f(y_{1}^m)).
\end{equation}

So by $f$ operation on both sides of \eqref{thgeneqn110} with
$f(x_1^m)$, we get,

\begin{equation} \label{thgeneqn2}
\begin{split}
f(f(\underbrace{\mathbf{0} , \ldots ,  \mathbf{0}}_{m-1} , f(x_1^m)),f(x_1^m),\underbrace{\mathbf{0} , \ldots ,  \mathbf{0}}_{m-2}) \preccurlyeq 
\\
f(f(\underbrace{\mathbf{0} , \ldots ,  \mathbf{0}}_{m-1} , f(x_1^m)) , f(y_1^m),\underbrace{\mathbf{0} , \ldots ,  \mathbf{0}}_{m-2}).
\end{split}
\end{equation}

So by $f$ operation on both sides of \eqref{thgeneqn110} with
$f(y_1^m)$, we get,

\begin{equation} \label{thgeneqn311}
\begin{split}
f(f(\underbrace{\mathbf{0} , \ldots ,  \mathbf{0}}_{m-1} , f(x_1^m)) , f(y_1^m),\underbrace{\mathbf{0} , \ldots ,  \mathbf{0}}_{m-2} ) \preccurlyeq 
\\
f(f(\underbrace{\mathbf{0} , \ldots ,  \mathbf{0}}_{m-1} , f(y_1^m)),\underbrace{\mathbf{0} , \ldots ,  \mathbf{0}}_{m-2},f(y_1^m)).
\end{split}
\end{equation}

From~\eqref{thgeneqn2} and \eqref{thgeneqn311}, we get,
\[f(\underbrace{\mathbf{0} , \ldots ,  \mathbf{0}}_{m-2} , \stackrel{(2)}{f(x_1^m)})  \preccurlyeq f(\underbrace{\mathbf{0} , \ldots ,  \mathbf{0}}_{m-2} , \stackrel{(2)}{f(y_1^m)}).\]

Similarly, we get for $m$ terms:
\[f(\stackrel{(m)}{f(x_1^m)})  \preccurlyeq f(\stackrel{(m)}{f(y_1^m)}).\]

\item[(ii)] We know that
\[z_j \preccurlyeq u_j,  \qquad 1 \leq j \leq n.\]

From Lemma \ref{lemma2} (ii), we get, 
\[g(z_1^n) \preccurlyeq g(u_1^n).\]
Which is represented as follows
\begin{equation} \label{thgeneqn410}
g(\underbrace{\mathbf{1}, \ldots , \mathbf{1}}_{n-1}, g(z_{1}^n)) \preccurlyeq g(\underbrace{\mathbf{1}, \ldots , \mathbf{1}}_{n-1}, g(u_{1}^n)).
\end{equation}

Now by $g$ operation on both sides of \eqref{thgeneqn410} with
$g(z_1^n)$, we get,

\begin{equation} \label{thgeneqn5}
g(\stackrel{(2)}{g(z_1^n)},\underbrace{\mathbf{1}, \ldots , \mathbf{1}}_{n-2}) \preccurlyeq g(g(z_1^n) , g(u_1^n),\underbrace{\mathbf{1}, \ldots , \mathbf{1}}_{n-2}).
\end{equation}

So by $g$ operation on both sides of \eqref{thgeneqn410} with
$g(u_1^n)$, we get,

\begin{equation} \label{thgeneqn6}
g(\underbrace{\mathbf{1}, \ldots , \mathbf{1}}_{n-2} , g(z_1^n) , g(u_1^n) ) \preccurlyeq g(\underbrace{\mathbf{1}, \ldots , \mathbf{1}}_{n-2} , \stackrel{(2)}{g(u_1^n)}).
\end{equation}

So now from \eqref{thgeneqn5} and \eqref{thgeneqn6}, we get,
\[g(\underbrace{\mathbf{1}, \ldots , \mathbf{1}}_{n-2} , \stackrel{(2)}{g(z_1^n)})  \preccurlyeq g(\underbrace{\mathbf{1}, \ldots , \mathbf{1}}_{n-2} , \stackrel{(2)}{f(u_1^n)}).\]

Similarly, we find for $n$ terms
\[g(\stackrel{(n)}{g(z_1^n)})  \preccurlyeq g(\stackrel{(n)}{g(u_1^n)}).\]

\item[(iii)]

From Lemma \ref{lemma2} (ii)
\[g(z_1^n) \preccurlyeq g(u_1^n).\]

Similar to part (i), we find $f$ operation of $m$ terms and get
\[f(\underbrace{g(z_1^n) , g(z_1^n) , \ldots  , g(z_1^n)}_m) \preccurlyeq f(\underbrace{g(u_1^n) , g(u_1^n) , \ldots , g(u_1^n)}_m),\]

\[f(\stackrel{(m)}{g(z_1^n)}) \preccurlyeq f(\stackrel{(m)}{g(u_1^n)}).\]

\item[(iv)] We know that
 \[x_i \preccurlyeq y_i ,  \qquad  1 \leq i \leq m,\]

 so from Lemma \ref{lemma2} (i), we get
\[f(x_1^m) \preccurlyeq f(y_1^m).\]

As proved in part (ii), we find $g$ operations of $n$ terms and get,
\[g(\underbrace{f(x_1^m) , f(x_1^m) , \ldots  , f(x_1^m)}_n) \preccurlyeq g(\underbrace{f(y_1^m) , f(y_1^m) , \ldots , f(y_1^m)}_n).\]
Thus, we get,
\[g(\stackrel{(n)}{f(x_1^m))} \preccurlyeq g(\stackrel{(n)}{f(y_1^m))}. \qedhere\]

\end{itemize}

\end{proof}

\begin{corollary} \label{corollary2}

  If $\preccurlyeq$ is a fault-tolerance partial order and $k < m$, $t
  < n$ where \(k , t \in \mathbb{Z}_+\), if \(x_i \preccurlyeq y_i\),
  \(z_j \preccurlyeq u_j\) for all disjoint components $x_i, z_j, y_i, u_j$, which are in $\mathcal{U}$, where $1 \leq i \leq m$ and $1 \leq j \leq n$ then:

\begin{itemize}

\item[(i)]
  \(f(\underbrace{\mathbf{0},\ldots,\mathbf{0}}_{m-k},\stackrel{(k)}{f(x_1^m)})
  \preccurlyeq
  f(\underbrace{\mathbf{0},\ldots,\mathbf{0}}_{m-k},\stackrel{(k)}{f(y_1^m)})\)
  
\item[(ii)]
  \(g(\underbrace{\mathbf{1},\ldots,\mathbf{1}}_{n-t}\stackrel{(t)}{g(z_1^n)})
  \preccurlyeq
  g(\underbrace{\mathbf{1},\ldots,\mathbf{1}}_{n-t},\stackrel{(t)}{g(u_1^n)})\).

\end{itemize}

\end{corollary}

\begin{proof}
\begin{itemize}

\item[(i)]
 Proof is similar to that of Theorem~\ref{thgen} (i), we find the $f$ operation of $k$ terms where \(\forall k \in \mathbb{Z}_+\), and $k < m$.
 \item[(ii)]
Proof is similar to that of Theorem~\ref{thgen} (ii), we find the $g$ operation of $t$ terms where \(\forall t \in \mathbb{Z}_+\), and $t < n$.  \qedhere
\end{itemize}
\end{proof}

We propose the following Theorem for very complex systems.
\begin{theorem} \label{thgen2}

  If $\preccurlyeq$ is a fault-tolerance partial order, disjoint
  components $a_{i}, b_{i}, c_{j} , d_{j}$, $x_{k}, y_{k}, z_{t}, u_{t}$
  are in $\mathcal{U}$ and $a_i \preccurlyeq b_i, c_j \preccurlyeq
  d_j$, $x_k \preccurlyeq y_k$ and $z_t \preccurlyeq u_t$, where $1
  \leq i \leq m, 1 \leq j \leq n, 1 \leq k \leq m$ and $1 \leq t \leq
  n$, then:

\begin{itemize}

\item[(i)] $f(x_1^{k-1} , f(a_{1}^m) , x_{k+1}^{m}) \preccurlyeq
  f(y_1^{k-1} , f(b_{1}^m) , y_{k+1}^{m})$,

for all $1 \leq k \leq m$; and 

\item[(ii)] $f(x_1^{k-1} , g(c_{1}^n) , x_{k+1}^{m}) \preccurlyeq
  f(y_1^{k-1} , g(d_{1}^n), y_{k+1}^{m})$,

for all $1 \leq k \leq m$; and

\item[(iii)] $g(z_1^{t-1} , f(a_{1}^m) , z_{t+1}^{n}) \preccurlyeq
  g(u_1^{t-1} , f(b_{1}^m) , u_{t+1}^{n})$,

for all $1 \leq t \leq n$; and 

\item[(iv)] $g(z_1^{t-1} , g(c_{1}^n) , z_{t+1}^{n}) \preccurlyeq
  g(u_1^{t-1} , g(d_{1}^n) , u_{t+1}^{n})$,

for all $1 \leq t \leq n$.

\end{itemize}

\end{theorem}

\begin{proof}

\begin{itemize}

\item[(i)]

  From Lemma \ref{lemma2} (i), if $a_i \preccurlyeq b_i$, then
  $f(a_1^m) \preccurlyeq f(b_1^m)$ for all $1 \leq i \leq m$.

  We prove in a similar manner as Lemma \ref{lemma2} (i) that
\[f(f(a_{1}^m),x_1,\underbrace{\mathbf{0} , \ldots ,  \mathbf{0}}_{m-2}) \preccurlyeq f(f(b_{1}^m) , y_1 , \underbrace{\mathbf{0} , \ldots ,  \mathbf{0}}_{m-2}).\] 

Similarly, we get,

\[f(f(a_{1}^m),x_{1}^{k-1},x_{k+1}^m) \preccurlyeq f(f(b_{1}^m) , y_{1}^{k-1} ,y_{k+1}^{m}).\] 
Thus,
\[f(x_{1}^{k-1},f(a_{1}^m),x_{k+1}^m) \preccurlyeq f(y_{1}^{k-1} , f(b_{1}^m) , y_{k+1}^{m}).\] 

Similar to the above, we can prove (ii), (iii) and (iv).  \qedhere

\end{itemize}

\end{proof}

\end{document}